\documentclass[reqno,12pt]{amsart}
\makeatletter
\renewcommand{\@biblabel}[1]{#1.}
\makeatother
\newcommand{\ssp}{\discretionary{}{}{\,}}         

\theoremstyle{plain}
\newtheorem{prp}{Proposition}
\newtheorem{lemma}[prp]{Lemma}

\newcommand{\BWM}[2][(r,q)]{\mathrm{BWM}{#1}_{#2}}
\newcommand{\qm}{\hat{q}}
\newcommand{\qnum}[1]{\lbrack \! \lbrack #1 \rbrack \! \rbrack }

\begin{document}
\title[]{Symmetrizer and Antisymmetrizer of the Birman--Wenzl--Murakami Algebras}
\author[]{Istv\'an Heckenberger and Axel Sch\"uler}
\thanks{Supported by the Deutsche Forschungsgemeinschaft}

\address{Universit\" at Leipzig, Mathematisches Institut,
 Augustusplatz 10, D-04109 Leipzig,
Germany}
\email{\{heckenbe,schueler\}@mathematik.uni-leipzig.de}

\begin{abstract}
Explicit formulas for the symmetrizer and the antisymmetrizer of the
Birman--Wenzl--Murakami algebras $\BWM n$ are given.
\end{abstract}

\maketitle


\section{Introduction}

The Birman--Wenzl--Murakami algebra was first defined and independently
studied by Birman and Wenzl \cite{a-BirWen89} and Murakami \cite{a-Murakami90}.
The Iwahori--Hecke algebras of Type A and the Birman--Wenzl--Murakami algebras 
naturally arise as centralizer algebras of tensor product corepresentations
of quantum groups of Type A and of Type B, C, and D, respectively
\cite{a-Wenzl88,a-Wenzl90}.
Irreducible characters and primitive idempotents of Hecke algebras have been
intensively studied during the last 30 years. 
There are explicit constructions of analogues of the Young symmetrizer
and Specht modules
\cite{a-Gyoja86, a-LeducRam}. 
However, in case of the Birman--Wenzl--Murakami algebras there are only
a few papers in this direction \cite{a-BirWen89,a-LeducRam}. 
Young symmetrizers for the Birman--Wenzl--Murakami algebra are not known.

The aim of this Letter is to give recursive but explicit formulas for the
most important minimal central idempotents
of the Birman--Wenzl--Murakami algebras
namely for the symmetrizer $S_n$
and the antisymmetrizer $A_n$.
For the symmetrizer,  we establish the formulas
$$
c S_n=S_{n-1}b^+_{n-1,1}=a^+_{n-1,1}S_{n-1}
=S_{n-1}[1]b^+_{1,n-1}=a^+_{n-1,1}S_{n-1}[1],
$$
where $c$ is a normalization constant.
Here $S_{n-1}[1]$ denotes the symmetrizer $S_{n-1}$ shifted by one.
The elements $a^+_{1,n-1}$ and $a^+_{n-1,1}$ are the substitutes for
the well-known sums of shuffle permutations 
$$
1 +s_1+s_1s_2+\cdots+ s_1\cdots
s_{n-1}
\quad
\text{ and }\quad 
1+s_{n-1}+s_{n-1}s_{n-2}+\cdots+s_{n-1}\cdots s_1
$$ in the group
algebra of the symmetric group.
The elements $b^+_{n-1,1}$ and $b^+_{1,n-1}$ correspond to the sum of the
inverse shuffle permutations.
Taking the quotient of
$\BWM n$ by the relations $e_1=\cdots  =e_{n-1}=0$ one gets the well-known
formula for the symmetrizer of the Hecke algebra
$$
S_n=q^{-\frac{n(n-1)}{2}}(\qnum{n}!)^{-1}\sum_{w} q^{\ell(w)}T_w,
$$
 where the
sum is over all elements $w$ of the symmetric group of $n$ variables.
Here $T_w$ denotes the corresponding elements of the Hecke algebra and
$\ell (w)$ is the length of $w$. Similar formulas
for the antisymmetrizer $A_n$ are also given.

\section{The Birman--Wenzl--Murakami Algebra}

Let $n\in \mathbb{N}$, $n\geq 2$, and $r,q\in \mathbb{C}\setminus \{0\}$.
Further, we assume that
$q^k\not= 1$ for any $k\in \mathbb{N}$ and
$r\not= \pm q^k$ for any $k\in \mathbb{Z}$.
We use the abbreviation $\qm $ for the complex number
$q-q^{-1}$. The symbol $\qnum{k}$, $k\in \mathbb{Z}$, denotes the complex
number $(q^k-q^{-k})/(q-q^{-1})$.

The Birman--Wenzl--Murakami algebra $\BWM n$ is the unital algebra 
over the complex numbers generated by the elements
$g_i,e_i$, $i=1,\ldots ,n-1$, and relations
\begin{align}\label{eq-bwmrel1}
g_ig_{i+1}g_i&=g_{i+1}g_ig_{i+1},& g_ig_j&=g_jg_i\quad \text{for $|i-j|>1$},\\
\label{eq-bwmrel2}
g_i^2&=1+\qm g_i -r^{-1}\qm e_i,& g_i e_i&=e_ig_i=r^{-1}e_i,\\
\label{eq-bwmrel3}
g_ig_{i+1}e_i&=e_{i+1}e_i, & \makebox[0pt][r]{$g_{i+1}g_ie_{i+1}$}&=e_ie_{i+1}.
\end{align}
For any $n\in \mathbb{N}$, $n\geq 2$, the algebra $\BWM n$ is naturally
embedded into the algebras $\BWM{n+k}$, $k\geq 0$.

There is an algebra automorphism $\alpha _n$ of $\BWM n$ defined by
$\alpha _n(g_i)=g_{n-i}$, $\alpha _n(e_i)=e_{n-i}$ for all $i=1,\ldots ,n-1$.
There is an algebra antiautomorphism $\beta _n$ of $\BWM n$ defined by
$\beta _n(g_i)=g_i$, $\beta _n(e_i)=e_i$ for any $i=1,\ldots ,n-1$.
There is an algebra isomorphism
$\gamma _n:\BWM n \to \BWM[(r,p)]{n}$, $p=-q^{-1}$, given by
$\gamma _n(g_i)=g_i$, $\gamma _n(e_i)=e_i$.

Let $s: \BWM{n}\to \BWM{n+1}$ be the algebra homomorphism defined by
$s(g_i):=g_{i+1}$, $s(e_i)=e_{i+1}$.
If $b=:b[0]\in \BWM n$ and $s^k(b)\in \BWM n$ for some $k\in \mathbb{N}$, then we
write $b[k]$ for the element $s^k(b)$.

\section{The Symmetrizer}

The symmetrizer $S_n$ is the unique nonzero element of the algebra
$\BWM n$ such that $S_n^2=S_n$ and $S_ng_i=qS_n$
for any $i=1,\ldots ,n-1$ \cite[Theorem 5.14]{a-LeducRam}.
It follows that $S_n$ is central, in particular $g_iS_n=qS_n$ for all
$i=1,\ldots ,n-1$.
It is well known that
\begin{align}\label{eq-S2}
S_2&=\frac{1}{q\qnum{2}}\left( 1+q g_1+\frac{q \qm }{1-qr}e_1\right).
\end{align}
Now let us introduce the elements
\begin{align}\label{eq-dk}
d^+_{k,i}&:=e_{k-1} e_{k-2} \cdots e_i
 \sum _{j=0}^{i-1} q^j g_{i-1}\cdots g_{i-j},
& 1\leq i<k,\\
b^+_{k,1}&:=\sum _{i=0}^k q^i g_k g_{k-1} \cdots g_{k+1-i}
+\frac{\qm }{1-q^{2k-1}r}\sum _{i=1}^k q^{2k-2i+1} d^+_{k+1,i},& 0\leq k<n.
\end{align}
More precisely, formula (\ref{eq-dk}) reads as $d^+_{k,i}=e_{k-1} e_{k-2} \cdots e_i \sum _{j=0}^{i-1} q^j \prod_{l=1}^j g_{i-l}$.
 
\begin{prp}
For any $n\geq 3$ the formula $S_n=S_{n-1}b^+_{n-1,1}/(q^{n-1}\qnum{n})$ holds.
\end{prp}

\begin{lemma}\label{l-dkgl}
The expression $S_{n-1}d^+_{n,k} g_l$, $l<n$, equals
\begin{align*}
& q S_{n-1}d^+_{n,k}-qr^{-1}\qm S_{n-1}e_{n-1} e_{n-2} \cdots e_l & &
\text{for $l<k$},\\
& qS_{n-1}d^+_{n,k-1}+r^{-1}S_{n-1}e_{n-1} e_{n-2} \cdots e_l & &
\text{for $l=k>1$},\\
& r^{-1} S_{n-1}d^+_{n,k} & &\text{for $l=k=1$},\\
& q^{-1}S_{n-1}d^+_{n,k+1}+\qm S_{n-1}d^+_{n,k}
-q^{2k-1}S_{n-1}e_{n-1}e_{n-2}\cdots e_l & &\text{for $l=k+1$},\\
& q S_{n-1}d^+_{n,k} & &\text{for $l\geq k+2$.}
\end{align*}
\end{lemma}

\begin{proof}[Proof of the Proposition]
Let $S'_n:=S_{n-1}b^+_{n-1,1}/(q^{n-1}\qnum{n})$.
We prove that $S'_ng_i=qS'_n$ for any $i=1,2,\ldots,n-1$. Then we conclude
that $S'_ne_i=q^{-1}S'_ng_i e_i=q^{-1}r^{-1}S'_ne_i$ and hence
$(q^{-1}r^{-1}-1)S'_ne_i=0$. Since $r\not= q^{-1}$, we obtain $S'_ne_i=0$.
Now (\ref{eq-S2}) and the assertion of the Proposition for $k<n$ imply
that
\begin{align}
{S'_n}^2=&S'_nS_2 b^+_{2,1}b^+_{3,1}\cdots b^+_{n-1,1}/(q^2\qnum{3}
q^3\qnum{4}\cdots q^{n-1}\qnum{n})=S'_n.
\end{align}
Hence, $S'_n$ is the symmetrizer of the algebra $\BWM n$.

Let $\tilde{S}_n=S_{n-1}b^+_{n-1,1}$. We prove that
$\tilde{S}_ng_i=q\tilde{S}_n$ for any $i=1,2,\ldots ,n-1$.
For this we use the relations of the algebra $\BWM n$ and Lemma \ref{l-dkgl}.
Observe that in Lemma \ref{l-dkgl} only $S_{n-1}$ occurs. Therefore, we
can use that $S_{n-1}g_i=qS_{n-1}$ for all $i=1,\ldots ,n-2$.
Let now $i\geq 2$. Then
\allowdisplaybreaks
\begin{align*}
\tilde{S}_ng_i=& S_{n-1}\left( \sum _{j=0}^{n-1}q^jg_{n-1}g_{n-2}\cdots g_{n-j}
+\frac{\qm }{1-q^{2n-3}r}\sum _{j=1}^{n-1}q^{2n-2j-1}d^+_{n,j}\right) g_i\\
=& S_{n-1}\left( q\sum _{j=0}^{n-i-2}q^jg_{n-1}\cdots g_{n-j}
+q^{n-i-1}g_{n-1}\cdots g_{i+1}g_i+\right. \\ \nopagebreak
&+q^{n-i}g_{n-1}\cdots g_{i+1}(1+\qm g_i -r^{-1}\qm e_i)
+q\sum _{j=n-i+1}^{n-1}q^jg_{n-1}\cdots g_{n-j}+\\ \nopagebreak
&+\frac{\qm }{1-q^{2n-3}r}\left( \sum _{j=1}^{i-2}q^{2n-2j-1}qd^+_{n,j}
+q^{2n-2i+1}(q^{-1}d^+_{n,i}+\qm d^+_{n,i-1}-\right. \\ \nopagebreak
&-q^{2i-3}e_{n-1}\cdots e_i)
+q^{2n-2i-1}(qd^+_{n,i-1}+r^{-1}e_{n-1}\cdots e_i)+\\ \nopagebreak
&\left. \left.
+\sum _{j=i+1}^{n-1}q^{2n-2j-1}(qd^+_{n,j}-qr^{-1}\qm e_{n-1}\cdots e_i)
\right) \right)\\
=& S_{n-1}\bigg( q\sum _{j=0}^{n-1}q^jg_{n-1}\cdots g_{n-j}
-q^{n-i}q^{-n+i+1}r^{-1}\qm e_{n-1}\cdots e_{i+1}e_i+ \\ \nobreak
&+\frac{\qm }{1-q^{2n-3}r}\bigg(\sum _{j=1}^{n-1}q^{2n-2j}d^+_{n,j}
+(-q^{2n-2}+q^{2n-2i-1}r^{-1}-\\ \nopagebreak
&-r^{-1}\qm q^{n-i}\qnum{n-1-i})e_{n-1}\cdots e_i\bigg) \bigg)\\
=&qS_{n-1}b^+_{n-1,1}+\qm S_{n-1}(-qr^{-1}
+(-q^{2n-2}+r^{-1}q)/(1-q^{2n-3}r))e_{n-1}\cdots e_i\\
=&q\tilde{S}_n.
\end{align*}
The case $i=1$ can be proven similarly.
\end{proof}

\begin{proof}[Proof  of the Lemma]
Let $l<k$. Using the relations (\ref{eq-bwmrel1})--(\ref{eq-bwmrel3}) we obtain
\begin{align*}
S_{n-1}d^+_{n,k}g_l=&
S_{n-1}e_{n-1}\cdots e_k\sum _{i=0}^{k-1} q^i g_{k-1}\cdots g_{k-i}g_l\\
=&S_{n-1}e_{n-1}\cdots e_k\left( q\sum _{i=0}^{k-l-2}q^ig_{k-1}\cdots g_{k-i}
+q^{k-l-1}g_{k-1}\cdots g_{l+1}g_l+\right. \\
&\left. +q^{k-l}g_{k-1}\cdots g_{l+1}(1+\qm g_l -r^{-1}\qm e_l)
+q\sum _{i=k-l+1}^{k-1}q^ig_{k-1}\cdots g_{k-i}\right)\\
=&S_{n-1}e_{n-1}\cdots e_k\left(
q\sum _{i=0}^{k-1}q^ig_{k-1}\cdots g_{k-i}
-q^{k-l}r^{-1}\qm q^{-k+l+1}e_{k-1}\cdots e_l\right)\\
=&qS_{n-1}d^+_{n,k}-qr^{-1}\qm S_{n-1}e_{n-1}e_{n-2}\cdots e_l.
\end{align*}
The other cases can be shown similarly.
\end{proof}

If we apply the automophism $\varphi :=\alpha _n$ (the antiautomorphisms
$\varphi :=\beta _n$ and $\varphi :=\alpha _n\circ \beta _n$, respectively,)
to the symmetrizer $S_n$, we reobtain $S_n$. Indeed, the image of $S_n$ is
an idempotent, $\varphi (S_n)\varphi (S_n)=\varphi (S_n^2)=\varphi (S_n)$.
Further,
\begin{gather}
\varphi (S_n)g_i=\varphi (S_ng_{n-i})=q\varphi (S_n)\qquad
\text{for all $i=1,\ldots ,n-1$}
\end{gather}
($\varphi (S_n)g_i=\varphi (\varphi (g_i)S_n)=\varphi (qS_n)$ for all
$i=1,\ldots ,n-1$).
Thus we obtain the following formulas for the symmetrizer $S_n$:
\begin{align}
{d'}^+_{k,i}:=&\alpha _k(d^+_{k,i})
=e_1e_2\cdots e_{k-i}\sum _{j=0}^{i-1}q^jg_{k+1-i}\cdots g_{k-i+j},\\
b^+_{1,k}:=&\alpha _{k+1}(b^+_{k,1})
=\sum _{i=0}^kq^ig_1g_2\cdots g_i
+\frac{\qm }{1-q^{2k-1}r}\sum _{i=1}^kq^{2k-2i+1}{d'}^+_{k+1,i},\\
\bar{d}^+_{k,i}:=&\beta _n(d^+_{k,i})=\sum _{j=0}^{i-1}q^j
g_{i-j}\cdots g_{i-1} e_ie_{i+1}\cdots e_{k-1},\\
a^+_{k,1}:=&\beta _n(b^+_{k,1})=\sum _{i=0}^k q^i g_{k+1-i}\cdots g_{k-1} g_k
+\frac{\qm }{1-q^{2k-1}r}\sum _{i=1}^k q^{2k-2i+1}\bar{d}^+_{k+1,i},\\
\bar{d}'{}^+_{k,i}:=&\beta _n({d'}^+_{k,i})
=\sum _{j=0}^{i-1} q^j g_{k-i+j}\cdots g_{k-i+1} e_{k-i}\cdots e_2 e_1,\\
a^+_{1,k}:=&\beta _n(b^+_{1,k})
=\sum _{i=0}^k q^i g_i\cdots g_2 g_1+\frac{\qm }{1-q^{2k-1}r}
\sum _{i=1}^k q^{2k-2i+1} \bar{d}'{}^+_{k+1,i}.
\end{align}
\begin{align}\notag
S_n&=S_{n-1}b^+_{n-1,1}/(q^{n-1}\qnum{n})
=S_{n-1}[1]b^+_{1,n-1}/(q^{n-1}\qnum{n})\\
&=a^+_{n-1,1}S_{n-1}/(q^{n-1}\qnum{n})
=a^+_{1,n-1}S_{n-1}[1]/(q^{n-1}\qnum{n}).
\end{align}

\section{The Antisymmetrizer}

It is well known that the antisymmetrizer $A_n$ is the unique
nonzero idempotent
in the algebra $\BWM n$ such that $A_n g_i=-q^{-1}g_i$ for all
$i=1,\ldots ,n-1$. Let us examine the image $A'_n$ of the symmetrizer
$S_n$ under the isomorphism $\gamma _n: \BWM[(r,p)]{n}\to \BWM n$,
$p=-q^{-1}$ \cite[Proposition 3.2 (c)]{a-Wenzl90}.
Obviously, ${A'_n}^2=\gamma _n(S_n)^2=\gamma _n(S_n^2)
=\gamma _n(S_n)=A'_n$, hence $A'_n$ is an idempotent. Moreover,
\begin{gather}
A'_n g_i=\gamma _n(S_n)\gamma _n(g_i)=\gamma _n(S_ng_i)=p\gamma _n(S_n)
=-q^{-1}A'_n
\end{gather}
for all $i=1,\ldots ,n-1$. The explicit formulas for $A'_n$ show that
$A'_n$ is nonzero. Hence, $A'_n$ is the antisymmetrizer of the algebra
$\BWM n$. We obtain the following formulas:
\allowdisplaybreaks
\begin{align}
d^-_{k,i}:=&\gamma _n(d^+_{k,i})=e_{k-1} e_{k-2} \cdots e_i
 \sum _{j=0}^{i-1} (-q)^{-j} g_{i-1}\cdots g_{i-j},\\
b^-_{k,1}:=&\gamma _n(b^+_{k,1})
=\sum _{i=0}^k (-q)^{-i} g_k g_{k-1} \cdots g_{k+1-i}
-\frac{\qm }{1+q^{-2k+1}r}\sum _{i=1}^k q^{2i-2k-1} d^-_{k+1,i},\\
{d'}^-_{k,i}:=&\gamma _n({d'}^+_{k,i})
=e_1e_2\cdots e_{k-i}\sum _{j=0}^{i-1}(-q)^{-j}g_{k+1-i}\cdots g_{k-i+j},\\
b^-_{1,k}:=&\gamma _n(b^+_{1,k})
=\sum _{i=0}^k(-q)^{-i}g_1g_2\cdots g_i
-\frac{\qm }{1+q^{-2k+1}r}\sum _{i=1}^kq^{2i-2k-1}{d'}^-_{k+1,i},\\
\bar{d}^-_{k,i}:=&\gamma _n(\bar{d}^+_{k,i})=\sum _{j=0}^{i-1}(-q)^{-j}
g_{i-j}\cdots g_{i-1} e_ie_{i+1}\cdots e_{k-1},\\
a^-_{k,1}:=&\gamma _n(a^+_{k,1})
=\sum _{i=0}^k (-q)^{-i} g_{k+1-i}\cdots g_{k-1} g_k
-\frac{\qm }{1+q^{-2k+1}r}\sum _{i=1}^k q^{2i-2k-1}\bar{d}^-_{k+1,i},\\
\bar{d}'{}^-_{k,i}:=&\gamma _n(\bar{d}'{}^+_{k,i})
=\sum _{j=0}^{i-1}(-q)^{-j} g_{k-i+j}\cdots g_{k-i+1} e_{k-i}\cdots e_2 e_1,\\
a^-_{1,k}:=&\gamma _n(a^+_{1,k})
=\sum _{i=0}^k (-q)^{-i} g_i\cdots g_2 g_1-\frac{\qm }{1+q^{-2k+1}r}
\sum _{i=1}^k q^{2i-2k-1} \bar{d}'{}^-_{k+1,i}.
\end{align}
For the antisymmetrizer $A_n=\gamma _n(S_n)$ the following equations hold:
\begin{align}
A_n&=q^{n-1}A_{n-1}b^-_{n-1,1}/\qnum{n}
=q^{n-1}A_{n-1}[1]b^-_{1,n-1}/\qnum{n}\\
&=q^{n-1}a^-_{n-1,1}A_{n-1}/\qnum{n}
=q^{n-1}a^-_{1,n-1}A_{n-1}[1]/\qnum{n}.
\end{align}

\providecommand{\bysame}{\leavevmode\hbox to3em{\hrulefill}\thinspace}

\end{document}